\documentclass[12pt]{amsart}
\usepackage{graphicx} 
\usepackage{amssymb} 
\usepackage{amsmath} 
\usepackage[a4paper, margin=1in]{geometry} 

\usepackage{amsthm} 
\theoremstyle{plain}
\newtheorem{proposition}{Proposition}[section]

\newtheorem{conjecture}[proposition]{Conjecture}

\newtheorem{definition}[proposition]{Definition}

\newtheorem{thm}[proposition]{Theorem}
\newtheorem{defn}[proposition]{Definition}

\usepackage{amsmath} 
\usepackage{amssymb} 
\usepackage{amsthm}
\usepackage{amsfonts}            
\usepackage{mathrsfs}
\usepackage[all]{xy}
\usepackage{mathtools}
\usepackage{mathabx}
\usepackage{colonequals}
\usepackage{tikz} 
\usepackage{microtype}

\usepackage{cite}
\usepackage{comment} 
\usepackage[pagebackref]{hyperref} 
\hypersetup{  
	colorlinks=true,
	linkcolor=red!50!black,,
	citecolor=green,
	urlcolor=magenta,
}

\usepackage{dynkin-diagrams}
\def\row#1/#2!{#1_{\IfStrEq{#2}{}{n}{#2}} & \dynkin{#1}{#2}\\}

\newcommand{\BA}{\mathbb A}

\newcommand{\BC}{\mathbb C}

\newcommand{\BE}{\mathbb E}
\newcommand{\BR}{\mathbb R}
\newcommand{\BX}{\mathbb X}
\newcommand{\BZ}{\mathbb Z}

\newcommand{\BQ}{\mathbb Q}
\newcommand{\BF}{\mathbb F}

\newcommand{\BV}{\mathbb V}
\newcommand{\BL}{\mathbb L}

\newcommand{\mfkd}{\mathfrak{d}}

\newcommand{\CM}{\mathcal M}
\newcommand{\CN}{\mathcal N}

\newcommand{\CE}{\mathcal E}

\newcommand{\CS}{\mathcal S}
\newcommand{\CL}{\mathcal L}
\newcommand{\CG}{\mathcal G}
\newcommand{\CZ}{\mathcal Z}

\newcommand{\U}{\mathrm U}

\newcommand{\GL}{\mathrm {GL}}
\newcommand{\SL}{\mathrm {SL}}

\newcommand{\Hom}{\mathrm {Hom}}
\newcommand{\End}{\mathrm {End}}
\newcommand{\pr}{\mathrm {pr}}
\newcommand{\Height}{\mathrm {Int}}
\newcommand{\Int}{\mathrm {Int}}

\newcommand{\Orb}{\mathrm {Orb}}

\newcommand{\Sh}{\mathrm {Sh}}

\newcommand{\Ch}{\mathrm {Ch}}
\newcommand{\Ker}{\mathrm {Ker}}
\newcommand{\Res}{\mathrm {Res}}

\newcommand{\Stab}{\mathrm {Stab}}
\newcommand{\Lie}{\mathrm {Lie}}

\newcommand{\Gal}{\mathrm {Gal}}
\newcommand{\Hk}{\mathrm {Hk}}
\newcommand{\rs}{\mathrm {rs}}

\newcommand{\Herm}{\mathrm {Herm}}
\newcommand{\Spec}{\mathrm {Spec}}
\newcommand{\Spf}{\mathrm {Spf}}

\title[TAFL, TAT and cycles on Shimura curves]{Non-reductive cycles and twisted arithmetic transfers for Shimura curves}
\author{Zhiyu Zhang}
\date{\today}

\begin{document}

\maketitle

\begin{abstract}
In this largely expository note, we explain some recent progress on new cycles on Shimura varieties and Rapoport-Zink spaces, (twisted) arithmetic fundamental lemma, and arithmetic analogs of relative Langlands program. We explain related formulations of arithmetic twisted Gan-Gross-Prasad conjecture, the proof of twisted AFL and certain arithmetic transfers.
\end{abstract}

\tableofcontents

\section{L-functions and cycles on Shimura varieties}

L-functions are fundamental objects in number theory. In relative Langlands program (e.g. \cite{BZSV}), conjecturally many automorphic L-functions could be expressed as \emph{period integrals} (``L=P'' patterns) up to \emph{explicit} local factors (assuming factorizations of global Arthur parameters on the Galois side). Waldspurger formula \cite{Waldspurger1985} and Ichino's triple product formula \cite{Ichino2008trilinear} are two basic examples related to inner forms of $\GL_2$. Gan-Gross-Prasad (GGP) conjectures \cite{GGP2011symplectic} are fundamental high-dimensional examples for tempered cuspidal automorphic forms with \emph{explicit} Ichino-Ikeda type refinements \cite{Ichino-Ikeda-orthogonal}, related local \emph{zeta integrals} and \emph{branching laws} (local GGP conjectures, e.g. \cite{He-GGP-discrete-series}) on real and $p$-adic groups, see \cite{gross2022road, BP2022ICM,Zhang2018ICM} for nice summaries. Theta liftings also provide (not necessarily tempered) fundamental examples, e.g. Rallis inner product formulas \cite{GQT2014-non-vanishing-global-Theta} and Saito–Kurokawa liftings \cite{Xue2018-Non-tempered-GGP}. In \cite{BZSV}, more examples are predicted for (graded) hyperspherical varieties, including strongly tempered spherical varieties.

There are many applications of these ``L=P'' patterns: detecting Langlands transfers (e.g. \cite{rallis1982langlands,rallis1990poles,Flicker-Rallis-base-change,FLO-unitary-period-stable-base-change,Ichino2022theta}); \emph{distribution}, \emph{non-vanishing}, \emph{algebraicity}, \emph{congruence}  of L-values; \emph{constructions} of $p$-adic L-functions.

Shimura varieties provide geometric tools to study arithmetic of (algebraic) automorphic forms and their $L$-functions. For example, (holomorphic) modular forms are special functions on the complex upper half plane $\mathbb H=(\mathbb C - \mathbb R)^+:=\{z \in \mathbb C | \mathrm{Im}(z) >0 \}$, and could be regarded as sections of certain line bundles on modular curves. In this way, their Ramanujan conjectures are proved via Weil conjecture \cite{Deligne-Ramanjuan-why-Shimura}, see also \cite{Brylinski-Labesse-1984-Asai-Hilbert,shahidi2011arthur-Ramanujan,clozel2013purity,Gan2023-Athur-Theta-Introduction}. Similar to ``L=P'' patterns, geometrically we have the fundamental \emph{Gross-Zagier formula} \cite{Gross-Zagier} for special cycles on modular curves. As a fundamental application, we know that an elliptic curve over $\mathbb Q$ has only finitely many $\mathbb Q$-points (resp. Mordell-Weil rank one) if its $L$-function is non-zero (resp. has order $1$) at the central point. We shall regarded speical cycles on Shimura varieties as geometric analogs of period integrals. From this perspective, we are interested in two topics:
 \begin{enumerate}
     \item \emph{Cycles on Shimura varieties}: construction and study of new cycles on Shimura varieties with arithmetic applications. We are interested in their \emph{distribution}, \emph{non-vanishing}, \emph{algebraicity}, \emph{modularity} (Kudla program) and \emph{intersection numbers}. There are applications to many areas, e.g. cohomology of (S-)arithmetic subgroups (e.g. \cite{millson1981geometric}), integral Euler systems (e.g. \cite{nekovavr1992values}), Noether-Lefschetz theory (e.g. \cite{Gromov-Witten-Noether-Lefschetz}), motivic Langlands transfers (e.g. \cite{Ichino-Prasana-JL-Archmiedean-Arthur}),  Hodge conjecture (e.g. \cite{blasius2000cohomology,BMM16}), Beilinson-Bloch-Kato conjecture (e.g. \cite{Li-Liu-proofs,LTXZZ,Proof-pAIPF-Disgeni-Liu,Proof-pAGGP-Disegni-Zhang}), and Tate conjecture over finite fields and number fields  (e.g. \cite{Xiao-Zhu-cycles,lemma2020algebraic}).
     \item \emph{Arithmetic relative Langlands program}: geometric analogs of ``L=P'' formulas relating L-functions to \emph{cycles on Shimura varieties}. There are recent ground-breaking progress for high-dimensional examples (with discrete series of minimal regular weights at infinity): \emph{arithmetic GGP conjecture}, \emph{arithmetic theta liftings} (Kudla program) and \emph{arithmetihc inner product formulas}, see \cite{Zhang-2019-incoherent-survey, 2024-AGGP-survey, 2024-ATheta-survey, Li-Liu-proofs} for nice summaries. Arthur's classification, local-global decompositions of height pairings, \emph{theta liftings and Weil representations}, \emph{explicit (analytic / mod $p$) geometry and cohomology (e.g. \cite[Lemma 7.3]{Li-Liu-proofs}\cite[Prop. 9.4.2]{Proof-pAGGP-Disegni-Zhang})) of Shimura varieties and Rapoport-Zink spaces} are crucial in the proofs. 
 \end{enumerate}

Firstly, to study Shimura varieties (over $\mathbb Q^{alg}$), we could use both complex (analytic) geometry and $p$-adic (analytic) geometry via base change along field embeddings
\[
\mathbb C \hookleftarrow \mathbb Q^{alg} \hookrightarrow \mathbb Q_p^{alg} \hookrightarrow \mathbb C_p.
\]
\begin{enumerate}
    \item \emph{Complex uniformization} studies Shimura varieties as groupoid quotient $\Sh_K(G,X_\infty)(\mathbb C)=G(\mathbb Q) \backslash X_\infty \times G(\mathbb A_f) /K$ with Hecke symmetry and local Galois symmetry at $\infty$ (complex Hodge theory), where differential forms and (sections of) metrized bundles naturally appear. Examples include Matsushima's formula, Falting's dual BGG complexes, and constructions of Green currents for special cycles \cite{Garcia2018superconnections,GarciaSankaran}.
    \item \emph{$p$-adic uniformization} (with the help of $p$-adic Hodge theory) gives new refined understandings of Shimura varieties with Hecke symmetry and local Galois symmetry at $p$, see e.g. \cite{Caraiani-Scholze-unitary,camargo2022-vanishing-completed,Igusa-p-adic-Matsushima, Fornea2024-plectic-uniformization}. We have a well-behaved theory of local Shimura varieties \cite{rapoport2014towards,Berkeley-Local-Sh} over $p$-adic numbers, with formal integral models called \emph{Rapoport-Zink spaces} \cite{RZ1996,pappas-rapoport-integralRZ}. For example, $\GL_2(\mathbb Q_p)$ acts naturally on the $p$-adic Drinfeld half plane $\mathbb H_p=\mathbb C_p - \mathbb Q_p$, and on its canonical integral model (a Rapoport-Zink space) over $\breve{\BZ}_p=W(\mathbb F_p^{alg})$, whose arithmetic quotients uniformize Shimura curves over $\breve{\BZ}_p$ \cite{2020-uniformization-Shimuracurve}. Moreover, the mod $p$ geometry of Shimura varieties (over $\mathbb F_p^{alg}$) gives more refined information (e.g. above Ramanujan conjecture, Langlands-Kottwitz method \cite{Kottwitz1990shimura}, and applications of function field Langlands and categorical local Langlands \cite{kisin2025strongly} \cite{2025genericcohomologyshimuravarieties}).
\end{enumerate}

Secondly, purely local methods are not sufficient to answer many global questions e.g. existence of automorphic forms and arithmetic of L-functions. To study \emph{global properties} of automorphic forms and Shimura varieties (with \emph{Hecke symmetry} and \emph{Galois symmetry}), the \emph{relative trace formula (RTF) method} is quite useful. Over number fields, trace formulas are powerful tools to prove the existence of transfers (e.g. cyclic base change \cite{arthur1989simple}, Jacquet-Langlands \cite{jacquet2006automorphic} and endoscopic classification \cite{Arthur2013}).

The idea of traces is fundamental for studying symmetry:
\[
\text{symmetry on spaces} \overset{\mathrm{Trace}}{\rightarrow} \text{(invariant) numbers and functions}. 
\]
For instance, characters are used to understand representations of finite groups. Classically, relative trace formulas for nice test functions are related to two period integrals (spectral expansion) we are interested in, and are decomposed into sums of products of (relative) \emph{local orbital integrals} (geometric expansion) over representatives of rational orbits (``relative Hitchin base''). Via trace formulas, global questions in Langlands program are related to purely local questions e.g. \emph{fundamental lemma} (FL) and \emph{transfers} (T) for \emph{comparison of orbital integrals}. For example, fundamental lemma for endoscopy proved by Ng\"{o} is crucial in establishing Arthur's classification of automorphic forms \cite{Arthur2013}\cite[\S 4]{Gan2023-Athur-Theta-Introduction}, and Langlands--Kottwitz method \cite{Kottwitz1990shimura} on quantitative descriptions of \emph{global Hecke-Galois actions} on Shimura varieties. Another example is the Jacquet-Rallis fundamental lemma proved in \cite{gordon2011fundamental,beuzart2021new} (whose arithmetic analog is the AFL discovered in \cite{Invent-AFL}), which is crucial for the recent proof of GGP conjecture for unitary groups. 

In this note, we explain how local and global methods are used together in arithmetic questions, at least for projective Shimura curves (see also \cite{KRY-Shimura-Curve,YZZ-Shimura-Curve,2020-uniformization-Shimuracurve}). 
\begin{enumerate}
    \item For topic (1): we construct new non-reductive cycles \cite{Proof-Invent-cycles-TAFL} e.g. global twisted Fourier-Jacobi cycles and local mirabolic cycles, which are related to geometric theta series. We formulate (ATGGP) \emph{arithmetic twisted GGP conjecture} on twisted \emph{Asai $L$-functions} \cite{TwistedGGP,Proof-Invent-cycles-TAFL}. For these new global cycles, we use that $H^1$ of $\U(n-1,1)$-Shimura varieties \cite{Liu-Fourier-Jacobi} contains Weil representations and theta liftings from $1$-dimensional skew-hermitian spaces (over finite adeles), building on \cite[Proposition 13.8]{BMM16} and Arthur's classification. Similar non-tempered cohomological results are known for $H^2$ of Siegel $3$-folds, $H^2$ of $\U(2,2) \sim SO(4,2)$ Shimura $4$-folds \cite{Ichino-Prasana-JL-Archmiedean-Arthur} and $H^4$ of Siegel $6$-folds \cite{sweeting2022tate}. ATGGP conjecture has arithmetic applications for Bloch-Kato conjecture on twisted Asai motives e.g. from (twisted) products of elliptic curves and CM abelian varieties, see \cite{Liu-TwistedTriple,Liu-Fourier-Jacobi,Proof-Invent-cycles-TAFL}.

    \item For topic (2): we review a powerful global method invented in \cite{Proof-AFL} for proving local arithmetic conjectures, e.g. arithmetic fundamental lemma (AFL) \cite{Proof-AFL-2}, and arithmetic transfers (AT) \cite{Proof-AT, LMZ-2025unitary}, where \emph{we bypass complicated questions on local geometry using automorphic representation theory (with explicit levels), and do induction based on simple known cases}. Here it's crucial to use slice versions of global automorphic functions and suitable (Gaussian) test functions at archimedean places to simplify the $\SL_2$-spectra. We explain a further generalization to prove \emph{twisted arithmetic fundamental lemma} (TAFL) involving general linear Rapoport--Zink spaces and local mirabolic cycles \cite{Proof-Invent-cycles-TAFL}, and certain twisted arithmetic transfers (TAT) for vertex hermitian lattices (Conjecture \ref{TAT(g,u) conj}), which are both crucial for RTF approach towards the ATGGP conjecture.
\end{enumerate}

In Lie theory, we have the fundamental classification of simple reductive groups / Lie algebras into ABCDEFG types (via simple roots / Dynkin diagrams). In this note, we just hope to emphasize the use of non-reductive (``type $P$'') groups (e.g. parabolic groups, Jacobi groups and mirabolic groups) and non-reductive cycles in arithmetic geometry and representation theory. Geometry of locally symmetric spaces beyond the case of Shimura varieties is also interesting, although we do not discuss it here. It shall be clear that there are many important topics and works of many people that are not mentioned in this note.

\subsection*{Acknowledgements}
I thank the organizers K. Namikawa and K. Gunji on the 2025 RIMS conference “Arithmetic aspects of automorphic forms
and automorphic representations” for their invitations and hospitality, and participants for sharing many related interesting mathematics. I also thank the organizers J.-B. Bost and S. Zhang on Simons Symposium on Periods and L-values of Motives (2024) for their invitations to give a related talk. I thank Y. Liu, W. Zhang for helpful comments.

\section{Producing cycles: Hecke translations, Weil representations and mirabolic cycles}

Asai motives occur naturally in the study of Hilbert modular varieties and unitary Shimura varieties.

Let $E_0/\mathbb Q$ be a real quadratic field, and $f$ be a Hilbert modular form over $E_0$ of weight $(2,2)$. By \cite{Brylinski-Labesse-1984-Asai-Hilbert}, the middle intersection cohomology of minimal compactification of Hilbert modular surfaces for $\GL_{2,E_0}$ (note that its reflex field is $\mathbb Q$) realizes the Asai motive $\mathrm{As}(f)=\mathrm{As}^+(f)$. Here $\mathrm{As}(f)$ is the (4-dimensional) $\Gal_{\mathbb Q}$-representation obtained as the \emph{tensor induction} of the $\Gal_{K}$-representation for $f$. The Asai (or twisted tensor) L-function $L(f, \mathrm{As}_{E_0/\mathbb Q}, s)=L(\mathrm{As}(f),s)$ is introduced in \cite{Asai1977} (see also \cite[\S 4]{Asai-classical-examples}). For applications to non-vanishing of Yoshida lifts, see \cite{Asai-Yoshida-Lift}. For Asai $L$-functions on $\GL_n$, see \cite{Flicker1988-Asai} for integral expressions, \cite{Asai-holomorphy-poles} for poles and \cite{BP-Asai-local-theory} for local L-factors. 

Similar to original Gross--Zagier formula, we study twisted Asai $L$-functions \cite{TwistedGGP} via base change along imaginary quadratic fields, and use unitary Shimura varieties. 

\subsection{Shimura varieties and cycles from Hecke translations}

More precisely, let $F/\mathbb Q$ (resp. $E_0/\mathbb Q$) be an imaginary (resp. real) quadratic field, and $E=F \otimes_\BQ E_0$ be the biquadratic extension of $F_0=\mathbb Q$. Let $V$ be a $n$-dimensional $F/\mathbb Q$ hermitian space of signature $(n-1,1)$ for a fixed embedding $\varphi_0: F \to \mathbb C$. Let $\mathbb D_{n-1} =\{ z \in \mathbb C^n | |z|<1\}=\U(n-1,1)(\mathbb R) / [ \U(n-1) \times \U(1)]$  be the complex unit ball of dimension $n-1$. Consider RSZ unitary Shimura varieties \cite{RSZ-AGGP}  for $V$ (resp.  $V_{E_0}$) over $E$ with suitable neat levels $K$ (resp. $K_G$)
\[
M(V) \to \Spec \,  E,  \quad M(V_{E_0}) \to \Spec \, E
\]
which are built from arithmetic quotients of $\mathbb D_{n-1}$ (resp. $\mathbb D_{n-1} \times \mathbb D_{n-1}$). For a discussion of reflex fields of unitary Shimura varieties, see \cite[(2.5)(2.10)(3.4), 4.3. Summary table]{RSZ-shimura}. We work over $E$ for simplicity, and $M(V),M(V_{E_0})$ may be defined over $F$. Choose compatible levels so the natural map $M(V) \to M(V_{E_0})$ is a closed embedding, which generalizes embeddings of Shimura curves to twisted Hilbert modular surfaces ($n=2$). 

\begin{defn}[Relative translation cycles]
Let $Z_1, Z_2$ be two closed subschemes in a scheme $X$. Given a correspondence $\pr_1 \times \pr_2: \Gamma \to X \times X$, let the \emph{relative $\Gamma$-translation} cycle of $Z_1$ in $Z_2$ be the fiber product
\[
\Gamma. Z_1 \cap_{X} Z_2= \{ x \in \Gamma | \pr_1(x) \in Z_1, \pr_2(x) \in Z_2 \} \to Z_1 \times Z_2.
\]  
If $Z_2=X$, we write $\pr_2: \Gamma. Z_1 \cap_{X} X \to X$ as the $\Gamma$-translation cycle of $Z_1$ in $X$. Define the {$\Gamma$-fixed point cycle} on $X$ as $\pr_1=\pr_2: \Gamma \cap_{X \times X} \Delta=\{ x \in \Gamma | \pr_1(x)=\pr_2(x) \} \to X$.
\end{defn}

Consider the Hecke correspondence $\Gamma_g=\Gamma_{K_GgK_G}$ on $X=M(V_{E_0})$ for $g \in \U(V_{E_0})(\mathbb A_{E_0,f})$. For $\varphi=\otimes_v \varphi_v \in \CS(\U(V_{E_0})(\mathbb A_{E_0,f}))^{K_G \times K_G}$, we have the Hecke correspondence $\Gamma_\varphi$ is a weighted sum of $\Gamma_g$ by $\varphi$. Define the \emph{twisted CM cycle}
\begin{equation} \label{TCM over E}
\pr_1: \mathrm{TCM}(\varphi):= \Gamma_\varphi. M(V) \cap_{M(V_{E_0})} M(V) \to M(V).    
\end{equation}
We have a natural decomposition (via moduli descriptions)
\[
\mathrm{TCM}(\varphi)=\coprod_{\alpha} \mathrm{TCM}_{\alpha}(\varphi),
\]
where $\alpha$ runs over $\U(V)(F_0) \backslash \U(V)(E_0) / \U(V)(F_0)$, or over $F_0$-points the GIT quotient $[\U(V) \backslash \U(V_{E_0})/\U(V)]$. If $E_0=\mathbb Q \times \mathbb Q$, we recover big CM cycles in \cite[\S 7.4]{Proof-AFL}. For general $E_0$, there is no description of $[\U(V) \backslash \U(V_{E_0})/\U(V)]$ (or $[\GL(V) \backslash \GL(V_{E_0})/ \GL(V)]$) via characteristic polynomials, instead we use the symmetrization map $\U(V_{E_0})/\U(V) \to \U(V_{E_0})$. For symmetric spaces over $\mathbb C$ (instead of $F_0$), we have relative Chevalley isomorphisms on double GIT quotients \cite{richardson1982orbits} via little Weyl groups, see also \cite[\S 2-4]{leslie2022stabilization}.

Over an archimedean place $v: E \to \mathbb C$, by complex uniformization we see that $\mathrm{TCM}(\varphi)(\mathbb C)$ is related to fixed points of $g_\infty \in \U(n-1,1)(\mathbb R)$ on $\mathbb D_{n-1}$ (which will be non-empty singleton for generic $g_\infty$). For a general embedding of Shimura varieties $\Sh_H=\Sh_{K_H}(H, X_{H, \infty}) \to \Sh_G=\Sh_K(G,X_\infty)$,  we can write Hecke correspondence $\Gamma_g$ as
$$
\pr_i: \Gamma_g(\mathbb C)= G(\mathbb Q) \backslash [X_\infty \times \Hk_g^\infty(\mathbb C)] \to \Sh_K(G,X)(\mathbb C), (x, g_1, g_2) \mapsto (x, g_i),
$$
where $\Hk_g^\infty(\mathbb C):=\{ (g_1, g_2) \in (G(\mathbb A_f)/K)^2 | g_1^{-1}g_2 \in KgK \}$. Hence we find
\[
\Gamma_g. \Sh_H \cap_{\Sh_G} \Sh_H = G(\mathbb Q) \backslash \{ (x, g_1, g_2) | \exists (y_i, h_i) \in X_{H,\infty} \times H(\mathbb A_f)/K_H,  (x, g_i) \in G(\mathbb Q).(y_i, h_i), i=1,2  \}.
\]
Therefore, there exists $\alpha \in G(\mathbb Q)$ such that $\alpha.(y_1, h_1)=(y_2,h_2)$ and $\alpha$ is well-defined in $H(\BQ)\backslash G(\BQ) / H(\BQ)$, giving above decompositions of $\mathrm{TCM}(\varphi)$ into $\mathrm{TCM}_{\alpha}(\varphi)$.

 We will only work with $\varphi$ and $\alpha$ such that $\mathrm{TCM}_{\alpha}(\varphi)(\mathbb C)$ is a \emph{non-empty finite set} of CM points. Over a $p$-adic place $v: E \to \mathbb C_p$ with $\varphi_v=1_{K_{G,v}}$, we will see later that by $p$-adic uniformization, $\mathrm{TCM}(\varphi)(\mathbb C_p)$ is related to $\mathbb C_p$-points of the \emph{derived twisted fixed cycles} appearing in the TAT conjecture \ref{TAT(g,u) conj} formulated below.

\subsection{Global cycles: geometric theta series and twisted Fourier-Jacobi cycles}

Fix a non-trivial additive character $\psi: \mathbb A_{\mathbb Q} \to \mathbb C$. For a $1$-dimensional skew-hermitian space $W$ over $F$ and a conjugate-sympletic character $\mu: \mathbb A_F^\times \to \mathbb C$, using $\psi$ and $\mu$ we produce the Weil representation $\omega=\omega(W, \mu)$ for $\U(V) \times \U(W)$. Here $V \otimes W$ is regarded as an (unpolarized) symplectic vector space (or Hamiltonian space as in \cite{BZSV}). To find geometric analogs of these Weil representations (at least over finite adeles), we need to understand low degree cohomology of (projective) Shimura varieties. In this subsection, we assume that $M(V)$ and $M(V_{E_0})$ are projective. Otherwise, we replace them with suitable toroidal compactifications.

By definition, $\mu|_{\mathbb A_{\mathbb Q}^\times}$ is the quadratic character $\eta$ associated to $F/\mathbb Q$ by class field theory. Assume that $\mu$ has weight $1$ (i.e. $\mu_\infty(z)=z/\bar{z}$ or $\bar{z}/z$). By works of Shimura, we have the CM abelian variety $A_\mu$ over $F$ with CM by the number field $M_\mu$ generated by values of $\mu \circ |-|^{1/2}_F$ over finite adeles.

\begin{thm}[Liu's formula \cite{Liu-Fourier-Jacobi}]
Assume that $n>2$. 
\begin{enumerate}
    \item As $\U(V)(\mathbb A_f)$-modules, the Betti cohomology (taking inverse limits along all levels)
    $$
    H^1(M(V)(\mathbb C), \mathbb C)= \oplus_{(W_f,\mu)}\omega_f(W_f, \mu),
    $$
 where $\omega_f(W_f, \mu)$ is the Weil representation for weight $1$ conjugate-sympletic Hecke character $\mu$ and $1$-dimensional skew-hermitian space $W_f$ over $\mathbb A_{F,f}$. 
    \item The Albanese of $M(V)$ is (up to isogeny) an explicit direct sum of $(A_\mu)_E$. 
\end{enumerate}
When $n=2$, we have a similar decomposition of the endoscopic part of $H^1$ and Albanese.
\end{thm}

The proof of Liu's formula (see also \cite[Lemma 3.15]{Li-Liu-proofs}) uses relations between theta liftings (which are well-behaved in the stable range \cite{li1997-stablerange}), Arthur's classification \cite{Arthur2013}\cite[Theorem 1.7.1.]{KMSW-Arthur-unitary}\cite[Theorem 12.10, Proposition 13.3.]{BMM16}, and poles of $L$-functions \cite[\S 10-13]{BMM16}\cite{GJS2009poles,GSS1997poles,wu2013irreducibility}. 
As we are interested in cohomology of degree $1$, a local key point of the proof is an explicit Vogan-Zuckerman's classification \cite{VoganZuckerman-cohomological}. For $\U(n-1,1)(\BR)$, irreducible unitary cohomological representations are classified as $\pi_{a,b}$ ($a+b \leq n-1, a,b \in \mathbb Z_{\geq 0}$) with minimal Hodge type $h^{a,b}(\pi_{a,b})=1$, see \cite[\S 3]{SimonShin-cohomology-U}\cite[\S 5]{Gerbelli-2023-Arthur-Cohomological}. When $n>1$, $\U(n-1,1)(\mathbb R)$ has a unique proper parabolic $P$ up to conjugacy and $\pi_{a,b}$ ($a+b<n-1$) is an explicit Langlands quotient from $P$ \cite[\S 3]{Eisenstein-parabolic-unitary}\cite[\S 4-5]{ClozelBergeron-2005}. For explicit classifications on other groups, see e.g. \cite[\S 9]{Kudla2022-remarks-cohomology} ($SO_0(n-2,2)$), \cite{Sarnobat2019functorial-Sp4} ($Sp_4(\mathbb R)$) and  \cite{Venkataramana2001some} ($\U(2,3)$).

For previous study of Albanese of unitary Shimura varieties, see \cite{Oda1981note,Murty1992albanese,Schoen2014arithmetic}. For cohomology of general Shimura varieties, see \cite[\S 9]{Arthur-1989-cohomology-Sh-parameter}\cite[\S 6]{1994zeta}. The non-vanishing of first Betti number \cite{Millson1976first} is special to $\U(n-1,1)$-Shimura varieties. 

Using Liu's formula, we now introduce more cycles. Given a morphism $\phi: M(V) \to (A_\mu)_E$, we say $\phi$ is of type $W_f$ if the induced map on complex Betti cohomology $H^1$ has image in $\omega_f(W_f,\mu)$. Then $\phi$ could be regarded as an element in $\omega_f(W_f, \mu)$. For above $\varphi$ and $\phi$ (of type $W_f$), the \emph{twisted Fourier-Jacobi cycle} for $(\varphi,\phi)$ is (certain cohomological trivial modification of) the map
\[
\mathrm{FJ}(\varphi,\phi): M(V) \overset{(\Gamma_\varphi, \phi)}\to M(V_{E_0}) \times (A_\mu)_E.
\]
Let $\mathrm{FJ}(W_f) \subseteq Ch^*(M(V_{E_0}) \times (A_\mu)_E )$ be the $\U(V_{E_0})(\mathbb A_{E_0,f})$-submodule generated by images of $\mathrm{FJ}(\varphi,\phi)$ in the $\mathbb C$-coefficient Chow group of $M(V_{E_0}) \times (A_\mu)_E $.  These cycles are related to Kudla's geometric theta series on $M(V)$ via certain doubling constructions (and cup products $H^1(M(V)(\BC)) \times H^1(M(V)(\BC)) \to H^2(M(V)(\BC))$) \cite{Liu-Fourier-Jacobi,Proof-Invent-cycles-TAFL}, which show their relations to local derived twisted CM cycles (\ref{twisted fixed cycles}) in TAFL conjecture \ref{TAT(g,u) conj} below via $p$-adic uniformization.

\subsection{Local cycles: mirabolic cycles on Rapoport-Zink spaces}

Besides Hecke symmetry on towers, Rapoport-Zink spaces have an extra symmetry from automorphisms of framing objects. Such extra symmetry is crucial in local Langlands program, and will be used in our formulation of TAFL. 

Let $p$ be a prime number, and we use local notations. Let $F=\mathbb Q_{p^2}$ be the unramified quadratic extension of $F_0=\mathbb Q_p$ with non-trivial Galois involution $x \rightarrow \bar{x}$. Let $V$ be a $n$-dimensional hermitian space over $F/F_0$. Let $L \subseteq V$ be a vertex lattice, i.e. $L \subseteq L^\vee \subseteq p^{-1}L$. Let $t=\dim_{\mathbb F_{p^2}} L^\vee / L \in [0,n]$ be the type of $L$. We have an isomorphism $ \breve{\BZ}_p \otimes O_F \cong O_F \times O_F$. Let $S$ be a scheme over $\breve{\BZ}_p$ where $p$ is nilpotent. A basic tuple (resp. $L$-polarized tuple) of signature $(1,1)$ is the data $(X, \iota)$ (resp. $(X, \iota, \lambda)$ where
\begin{enumerate}
    \item $X$ is a $p$-divisible group over $S$ of height $2n$ and dimension $n$.
    \item $\iota: O_F \to \End(X)$ is an action s.t. $\Lie X$ has signature $(n-1,1)$ as $O_S \otimes O_F $-modules.
    \item (for $L$-polarized tuples) $\lambda: X \to X^\vee$ is a semi-linear $L$-polarization and $\Ker \lambda \in X[p]$ has order $p^{2t}$.
\end{enumerate}

Here polarizations are defined as anti-symmetric isogenies
i.e. $\lambda^\vee=-\lambda$, $\lambda \circ \iota(a)=\iota(\bar{a}) \circ \lambda (a \in O_F)$, see \cite[Page 9]{RSZ-Annalen} (called quasi-polarizations in \cite{oort2004foliations}). Fix a $L$-polarized basic tuple  $(\BX, \iota_{\BX}, \lambda_{\BX})$ over $\mathbb F_p^{alg}$ (framing tuple). Fix an $O_F$-polarized basic tuple $(\BE, \iota_\BE, \lambda_{\BE})$ of signature $(1, 0)$ over $\mathbb F_p^{alg}$ with canonical lifting $\CE$ to $\Spf \breve{\BZ}_p$. Consider the hermitian space of special quasi-homomorphisms
\[
\BV=\Hom_{O_F}(\BE, \BX) \otimes \BQ
\]
which has dimension $n$ over $F$. Note that $\BV \not \cong V$. Consider the unitary Rapoport-Zink space for $L$
\[
\CN_L \to \Spf \breve{\BZ}_p
\]   
as the moduli functor sending a test scheme $S$ where $p$ is nilpotent to isomorphism classes of $(X, \iota, \lambda, \rho)$ where $(X, \iota, \lambda)$ is a $L$-polarized basic tuple over $S$ and $\rho: X \times S/p \to \mathbb X \times S/p$ is a height $0$ $O_F$-linear quasi-isogeny compatible with polarizations.

Similarly, define the general linear Rapoport-Zink space
\[
\CN^\GL_n \to \Spf \breve{\BZ}_p
\]
as the moduli functor of $(X, \iota, \rho)$ where $(X, \iota)$ is a basic tuple over $S$ and $\rho: X \times S/p \to \mathbb X \times S/p$ is a height $0$ $O_F$-linear quasi-isogeny. We refer to the book \cite{RZ1996} on PEL type Rapoport-Zink spaces as deformation spaces of quasi-isogenies of $p$-divisible groups. By forgetting the polarization $\lambda$, we have a natural embedding
\[
\CN_L \to \CN_n^\GL
\]
is equivariant under the natural action of $\U(\BV)$ (resp. $\GL(\BV)$) on $\CN_L$ (resp. $\CN_n^\GL$)

Consider a non-zero vector $u \in \BV - 0$. The \emph{Kudla--Rapoport cycle} \cite{KR2011-local-cycle} \cite{Proof-AT}
$$
\CZ(u) \to \CN_L
$$ is the closed formal subscheme of $\CN_L$ sending a $\Spf O_{\breve{\BZ}_p}$-scheme $S$ to the subset $(X, \iota_X, \lambda_X, \rho_X) \in \CN_L(S)$ such that $(\rho_{X})^{-1}  \circ u: \BE \times_{\BF}  S/p \to X \times S/p $ 
	lifts to a homomorphism from $\CE$ to $X$ over $S$. By \cite[Prop 3.5]{KR2011-local-cycle}, $\CZ(u)$ is a Cartier divisor in $\CN_L$.

\begin{definition}
  The \emph{mirabolic special cycle} \cite{Proof-Invent-cycles-TAFL} (which is a Cartier divisor on $\CN^\GL_n$)
	$$
	\CZ^{\GL}(u) \to \CN^\GL_n
	$$ is the subfunctor of $\CN^\GL_n$ sending a $\Spf O_{\breve{\BZ}_p}$-scheme $S$ to $(X, \iota_X, \rho_X) \in \CN^\GL_n(S)$ such that
	\[ (\rho_{X})^{-1}  \circ u: \BE \times_{\BF} S/p \to X \times_{S}  S/p  \]  
	lifts to a homomorphism from $\CE$ to $X$ over $S$. For $u^* \in \BV^* -0 $,  the mirabolic special cycle
	$$
\CZ^{\GL^*}(u^*) \to \CN^\GL_n
$$ is the subfunctor of $\CN^\GL_n$ sending a $\Spf O_{\breve{\BZ}_p}$-scheme $S$ to $(X, \iota_X, \rho_X) \in \CN^\GL_n(S)$ such that $u^* \circ \rho_{X} :  X \times_{S} S/p \to \BE \times_{\BF} S/p$  lifts to a homomorphism from $X$ to $\CE$ over $S$.
  
\end{definition}
Via pullbacks from $\CN_n^\GL$, we obtain mirabolic special cycles on many other Rapoport-Zink spaces. For example, $\CZ^{\GL}(u)|_{\CN_L}=\CZ(u)$. The name of mirabolic cycles comes from that $\CZ^{\GL}(u)$ is stable under the action of the mirabolic subgroup $\Stab(u) \leq \GL(\mathbb V)$ and is an analog of symmetric spaces for mirabolic subgroups (over the generic fiber and the reduced locus). The fundamental work \cite{howard2019linear} shows the linear invariance of Kudla--Rapoport cycles (via Koszul complexs and deformation theory), and it is possible to prove similar things for mirabolic special cycles.

\section{L-functions and ATGGP conjecture via twisted Fourier-Jacobi cycles}

Let $\pi$ be an irreducible cuspidal automorphic representation of $\U(V_{E_0})(\mathbb A_{E_0})$. Let $\Pi$ be a (weak) base change of $\pi$ to $\GL_n(\mathbb A_{E})$. Let $\omega=\omega(\mu, W)$ be the Weil representation for $V \otimes W$. We formulate the (unrefined) conjecture on L-functions as follows.

\begin{conjecture}[Twisted GGP conjecture \cite{TwistedGGP}]
The following are equivalent:
\begin{enumerate}
    \item The period integral $\int_{[\U(V)]} \varphi(g) \overline{\phi(g)} \not =0 $ for some $\varphi \in \pi$ and $\phi \in \omega$.
    \item The twisted Asai L-value $L(1/2, \Pi, As_{E/F} \otimes \mu^{-1}) \not =0$ and $\Hom_{\U(V)(\mathbb A)}(\varphi, \omega) \not =0$.
\end{enumerate}
\end{conjecture}

\begin{conjecture}[Arithmetic twisted GGP conjecture \cite{Proof-Invent-cycles-TAFL}] 
Assume that $\pi$ is tempered and $\Pi$ is cuspidal. Moreover, assume that $\Pi$ is relevant, i.e.  for every archimedean place $v$ of $E$, $\Pi_v$ is isomorphic to the (irreducible) principal series representation induced by the characters $(arg^{1-n}, arg^{3-n}, ... , arg^{n-3},arg^{n-1})$ where $arg(z)=(\frac{z}{\bar{z}})^{1/2}: \mathbb C^\times \to \mathbb C^\times$. Then the following are equivalent:
\begin{enumerate}
\item $\Hom_{\U(V)(\mathbb A_{E_0,f})}(\pi_f, \mathrm{FJ}(W_f)) \not =0$.
\item $L'(1/2, \Pi, As_{E/F} \otimes \mu^{-1}) \not =0$ and $\Hom_{\U(V)(\mathbb A_f)}(\pi_f, \omega(\mu, W_f)) \not =0$.
\end{enumerate}
\end{conjecture}

See also the arithmetic Gan-Gross-Prasad conjecture formulated in \cite[Conjecture 6.10]{RSZ-AGGP} using cohomological tempered Arthur parameters. To formulate a refined version of ATGGP conjecture (using conjectural Beilinson-Bloch height pairings, see also recent progress \cite{zhang2021height}), see \cite{Liu-Fourier-Jacobi,Proof-Invent-cycles-TAFL}.  Given a good theory of twisted $p$-adic Asai L-functions, it is possible to formulate $p$-adic version of twisted Gan-Gross-Prasad conjecture and arithmetic analogs (see e.g. \cite{Proof-pAGGP-Disegni-Zhang}), where TAFL and TAT are still crucial.

As we only have Hecke actions over finite adeles on Shimura varieties, to formulate comparisons of arithmetic trace formulas (e.g. from Arakelov intersection pairings) and analytic trace formulas (e.g. from automorphic kernel functions), we often choose ($K_\infty$-finite) archimedean test functions by design, e.g. pseudo-coefficients \cite{ClozelDelorme-pseudocoeff} in the Langlands-Kottwitz method, and Gaussian test functions in the RTF approach towards AGGP conjecture \cite{Proof-AFL}.

In the case $\dim_{F} V=2$, ATGGP conjecture and twisted Fourier-Jacobi cycles are related to arithmetic triple product formulas, and arithmetic of Asai motives of elliptic curves over $E_0$ (e.g. from symmetric square of elliptic curves over $\BQ$ via base change). For similar examples, see \cite[\S 1.3]{Liu-TwistedTriple} (via orthogonal GGP conjecture) which does not use imaginary quadratic field and relies on certain Tate conjectures over finite fields \cite[Proposition 4.9]{Liu-TwistedTriple}, and also \cite[\S 1.5]{Liu-Fourier-Jacobi} (via untwisted GGP conjecture).

\section{Proof of TAFL: holomorphic Fourier expansion and local-global relations}

\subsection{Twisted arithmetic fundamental lemma (TAFL) and arithmetic transfers (TAT): formulations}

The arithmetic fundamental lemma (AFL) discovered by Wei Zhang \cite{Invent-AFL} (see also \cite{Liu-Fourier-Jacobi} for the Fourier-Jacobi case) and the subsequent arithmetic transfer (AT) conjectures \cite{RSZ-Duke,RSZ-Annalen,Proof-AT,CRZ-HeckeAFL,CRZ-quasiAFL} are identities between derived orbital integrals of certain test functions and arithmetic intersection numbers of diagonal cycles on unitary Rapoport-Zink spaces with certain levels. Proved in \cite{Proof-AFL,Proof-AFL-2,Proof-AT}, AFL and certain AT are used crucially in a recent proof of (refined) $p$-adic AGGP conjecture \cite{Proof-pAGGP-Disegni-Zhang} for a good ordinary split prime $p$, with applications to $p$-adic Beilinson-Bloch-Kato conjecture. Anticipating similar applications to arithmetic questions, we now give the formulation of twisted arithmetic fundamental lemma (TAFL) and certain twisted arithmetic transfers (TAT). 

Recall we have a $\U(\BV)$-equivariant closed embedding of regular formal schemes
\[
\CN_L \to \CN_n^\GL. 
\]
Recall for a reductive group action on an affine variety over a field, an orbit is called semi-simple (resp. regular) if the orbit is closed (resp. has trivial stabilizers). 

\begin{definition}\label{twisted fixed cycles}
For regular semi-simple $g \in \GL(\BV)(F_0)$, the derived twisted fixed cycle of $g$ is the derived cycle (in the Grothendieck group of $\CN_L$ with support in the schematic intersection $g \CN_L \cap_{\CN^\GL_n} \CN_L$)
\begin{equation}
\CN_L^{\Herm}(g)= g \CN_L \cap^{\BL}_{\CN^\GL_n} \CN_L \in K_0^{g \CN_L \cap_{\CN^\GL_n} \CN_L}(\CN_L),
\end{equation}    
which is a derived $1$-cycle by \cite[Lemma B.2]{Proof-AFL} and only depends on $g \in \GL(\BV)/\U(\BV)$. 
\end{definition}

Here for a Noetherian
formal scheme $X$ (resp. with a closed formal scheme $Y \subseteq X$), we write $K'_0(X)$ (resp. $K_0^{Y}(X)$) as the Grothendieck group of coherent sheaves of $O_X$-modules (resp. finite complexes of coherent locally free $O_X$-modules acyclic outside $Y$).

We identify $X_\BV:=\GL(\BV)/\U(\BV)=\Herm(\BV, h_\BV)$ as the set of hermitian structures $A$ on $\BV$ such that $(\BV, A) \cong (\BV, h_\BV)$. For any regular semi-simple pair $(g,u) \in (X_\BV \times \BV)(F_0)_\rs$, consider the derived intersection number
\[
\Height^{\Herm, \BV}(g, u) = \CN_L \cap^\BL_{\CN^\GL_n} g \CN_L \cap^\BL_{\CN^\GL_n} \CZ^\GL(u) \in \mathbb Q.
\]
Here we use that the schematic support is a proper scheme to define the intersection number via Euler characteristics of coherent sheaves.

Choose an orthogonal basis of $L \subseteq L^\vee$, then the basis generates a chain of $\mathbb Z_p$-lattices $L_0=\mathbb Z_p^n \subseteq L_0^\vee=p^{-1}\mathbb Z_p^t \oplus \mathbb Z_p^{n-t}$ in $\mathbb Q_p^n$. The stabilizer $\GL(L_0,L_0^\vee)$ is a parahoric subgroup of $\GL_n(\mathbb Q_p)$. The space  $\GL_n(\mathbb Q_p) \times \mathbb Q_p^n \times (\mathbb Q_p^n)^*$ has a natural action of $h \in \GL_n(\mathbb Q_p)$ by 
\[
h. (\gamma, u_1, u_2) = (h^{-1}\gamma h, h^{-1}u_1, u_2h).
\]
For regular semi-simple $(\gamma, u_1, u_2)$ and $(F', \phi') \in \CS(\GL_n(\mathbb Q_p) \times \mathbb Q_p^n \times (\mathbb Q_p^n)^*)$, we form orbital integrals ($s \in \mathbb C$) which converge absolutely 
\begin{equation}\label{defn: orb for TAT with transfer factor}
\Orb((\gamma, u_1,u_2), (F',\phi'),s):= \omega(\gamma, u_1, u_2) \int_{ \GL_n(\mathbb Q_p)} F'(h^{-1}\gamma h) \phi'(h^{-1}u_1,u_2h)
(-1)^{val(h)}|h|^s dh    
\end{equation}
with Haar measure on $\GL_n(\mathbb Q_p)$ normalized such that $\GL(L_0,L_0^\vee)$ has volume $1$ and $\omega(\gamma, u_1, u_2) \in \{\pm 1\}$ is a transfer factor. We form twisted arithmetic transfer (TAT) conjecture for $L$ as follows, see \cite{Proof-Invent-cycles-TAFL} (which is called TAFL when $t=0$).

\begin{conjecture}\label{TAT(g,u) conj}
($\mathrm{TAT}(g,u)$)
For any regular semi-simple $(g,u) \in (\GL(\BV)/\U(\BV) \times \BV)(\mathbb Q_p)_\rs$ matching $(\gamma,u_1,u_2) \in (\GL_n(\mathbb Q_p) \times \mathbb Q_p^n \times (\mathbb Q_p^n)^*)_\rs$, we have
\[			
\Height^{\Herm, \BV}(g,u) \log p = -  \frac{d}{ds}|_{s=0} \Orb( (\gamma,u_1, u_2), 1_{\GL(L_0,L_0^\vee)} \times 1_{L_0} \times 1_{(L_0^\vee)^*},s) \in \mathbb Q \log p.
\]
\end{conjecture}

The formulation is very similar to the global intersection problem of twisted  CM cycles and geometric theta series on unitary Shimura varieties, which will be the starting point of a global proof. 

Let's briefly explain the proof of TAFL \cite{Proof-Invent-cycles-TAFL} and TAT using the double induction method \cite{Proof-AT} in the case $n=2$ (to avoid non-properness of unitary Shimura varieties for $F_0=\mathbb Q$), i.e. we are working with Shimura curves. To obtain a proof for general $n$, we can work with projective RSZ Shimura varieties for totally real field $F_0 \not =\mathbb Q$, see \cite{Proof-AFL-2}\cite{Proof-Invent-cycles-TAFL}.

\subsection{Step I: globalization of intersections of cycles and supports}
Our proof is based on the idea of producing and comparing two \emph{holomorphic} Hilbert modular forms via comparing totally positive Fourier coefficients of finite parts \cite[Lemma 13.6]{Proof-AFL}. 

Fix a local regular semi-simple pair $(g_p, u_p)$ in Conjecture \ref{TAT(g,u) conj}. The case $n=1$ is true by direct computation. We do induction on $n$ and assume that Conjecture \ref{TAT(g,u) conj} is true in rank $n-1$. Let $F_0=\mathbb Q$. Choose $F=\mathbb Q(\sqrt{-D_1})$ ($D_1>0$), $E_0=\mathbb Q(\sqrt{D_2})$ $(D_2>0)$, and $E=\mathbb Q(\sqrt{-D_1},\sqrt{D_2})$, such that $p$ is inert in $E$ and $F$. Choose a $n$-dimensional $F/F_0$-hermitian space $V$ with a hermitian lattice $L$ such that 
\begin{enumerate}
    \item $L_p$ is self-dual.
    \item $V$ is of signature $(n-1, 1)$ for an embedding $\varphi_0: F \hookrightarrow \mathbb C$.
    \item There exists an inert place $w'$ of $F/F_0$ such that $V_{w'}$ is not split.
\end{enumerate}

Let $V^{(v_0)}=V^{(p)}$ be the totally definite nearby $F$-hermitian space of $V$ at the place $v_0=p$ of $F_0=\mathbb Q$. Then $V^{(v_0)} \otimes \BQ_p = \mathbb V_p$. 

By local constancy of orbital integrals and intersection numbers, we may choose a global regular semi-simple pair $(g_0,u_0) \in (\U(V^{(v_0)}) \times V^{(v_0)})(F_0)$ that is $p$-adic closely enough to $(g_p,u_p)$, such that 
\begin{enumerate}
    \item Conjecture \ref{TAT(g,u) conj} for $(g_p, u_p)$ and $(g_0,u_0)$ are equivalent.
    \item $\xi_0=(u_0,u_0) \in \mathbb Q_{>0}$ is totally positive. 
\end{enumerate}

We choose a large enough finite set $\mfkd$ of finite places of $F_0$ such that
\begin{itemize}
    \item $p$-adic places are not in $\mfkd$.
    \item any $w \not \in \mfkd$ is unramified in $F$ and $E_0$.
    \item $L_w$ is self-dual, for any $\ell$-adic place $w \not \in \mfkd$ ($\ell \not =p$) of $F_0$ inert in $F$.
    \item For any $\ell$-adic place $w \not \in \mfkd$ ($\ell \not = p$), the image $\alpha=\alpha(g_0)$ of $g_0$ in the GIT quotient $[\U(V^{(v_0)}_w)\backslash \U(V^{(v_0)}_{E_0,w}) / \U(V^{(v_0)}_w)](F_{0,w})$ is of maximal order, which is well-formulated as any $w \not \in \mfkd$ is unramified in $F$ and $E_0$.
\end{itemize}

We consider the RSZ integral model \cite{RSZ-AGGP,RSZ-AGGP}\cite[Appendix C.]{Liu-Fourier-Jacobi}
$$
\CM \to \Spec \, O_{F}[\mfkd^{-1}]
$$ 
of level $K$, where $K_w=\U(L_w)$ at any place $w \not \in \mfkd$, which is a $n$-dimensional \emph{regular} (from local models) and \emph{projective} scheme with generic fiber $M(V)$. Here we work with sufficient small levels at $\mfkd$. 

Via Serre tensor construction, we have a closed embedding of RSZ integral models for $M(V) \to M(V^{E_0})$:
\[
\CM \to \CM^{E_0}.
\]
with compatible level $K$ and $K_{G}$, such that $K_{G,w}=\U(L \otimes_{O_{F_0}} O_{E_0,w})$ for any place $w \not \in \mfkd$.

Let $\Ch^{1,adm}(\CM)$ be the subgroup of first arithmetic Chow group (i.e. arithmetic Picard group) of $\CM$ with admissible Green functions $\mathcal{G}$ at infinity, i.e. $\mathcal{G}$ has harmonic curvature with respect to the naturally metrized Hodge bundle $\omega_\mathbb C$ on $M(V)(\mathbb C)$. Let $\CZ_1(\CM)$ be the quotient of the group of $1$-cycles on $\CM$ by the subgroup generated by $1$-cycles
on $\CM$ that are contained in a closed fiber and rationally trivial within that fiber. We have a truncated Arakelov intersection pairing 
\[
(-,-)_{Ara}: \Ch^{1,adm}(\CM) \times \CZ_1(\CM) \to \mathbb R_\mfkd.
\]
Here the quotient $\BQ$-vector space $\BR_{\mfkd}=   \BR / \text{span}_{\BQ} \{ \log \ell | \exists w  \in \mfkd, \, w |\ell \}. $
A key feature is that if the support of intersection is empty on the generic fiber of $\CM$, then $(-,-)_{Ara}$ could be decomposed into local terms (hence related to our local intersection problems via uniformization) indexed by places $w \not \in \Delta$ of $F_0$:
\[
(-,-)_{Ara}=(-,-)_{Ara,\infty} + \sum_{w \nmid \infty} (-,-)_{Ara,w}. 
\] 

Consider a $\mathbb Q$-valued function $\varphi=1_{K_G^\mfkd} \times \varphi_\mfkd \in \CS(\U(V)(\mathbb A_{E,f}))^{K_G \times K_G}$. Then $\varphi_\mfkd$ gives a finite sum of prime-to-$p$ Hecke correspondences on $M(V_{E_0})$ which can be extended to the integral model (via moduli descriptions):
\[
\Gamma_{\varphi_\mfkd} \to \CM^{E_0} \times \CM^{E_0}.
\]
We form the \emph{derived twisted CM cycle} for $\varphi$ (viewed as an element in the Grothendieck group of coherent sheaves in the schematic intersections)
\[
{}^{\mathbb L} \mathcal{TCM}_{\alpha}(\varphi):= (\Gamma_{\varphi_\mfkd}. \mathcal{M} \cap^{\BL}_{\mathcal{M}^{E_0}} \mathcal{M} ) \to \mathcal{M}
\]
which on the generic fiber recovers $\mathrm{TCM}(\varphi)$. Denote by ${}^{\mathbb L} \mathcal{TCM}_{\alpha}(\varphi)$ the $\alpha$-part of ${}^{\mathbb L} \mathcal{TCM}(\varphi)$. The schematic intersection $\mathcal{TCM}_{\alpha}(\varphi) \to \mathcal{M}$ is finite and unramified (but may have large dimension components), see the proof of \cite[Proposition 7.9]{Proof-AFL}. Note that $\mathcal{TCM}(g_\mfkd)$ is of finite type hence it has only finitely many connected components, and $\mathcal{TCM}_\alpha(g_\mfkd)$ is an open and closed subscheme, see e.g. \cite[Lemma 7.11, 7.15]{Proof-AFL}. In terms of abelian schemes (ignoring the toric $Z^\BQ$-part), $\mathcal{TCM}(g_\mfkd)$ will roughly classify a pair $(A_1, A_2)$ of abelian schemes over a test scheme $S$ with polarizations, compatible $O_F$-actions of Kottwitz signature $(n-1,1)$ and level structures at $\mfkd$, with an $O_E$-linear quasi-isogeny $\varphi: A_1 \otimes_{O_F} O_E \to  A_2 \otimes_{O_F} O_E$ whose realization on $\mfkd$-adic rational Tate modules is compatible with $g_\mfkd \in \U(V_{E_0,\mfkd})$ under level structures.

Similarly, to globalize $\CZ(u_p)$ we consider the geometric theta series of Kudla-Rapoport divisors\cite{KR2011-local-cycle,KudlaRapoport-cycle-global}) for $\mathbb Q$-valued function $\phi \in \CS(V(\mathbb A_{f}))^{K}$:
\[
Z(h, \phi) \to M(V) , \quad h \in \SL_2(\mathbb A_{F_0}),
\]
which is a holomorphic modular form of weight $n$ with coefficients in $\mathrm{Ch}^1(M(V))$ by \cite{Liu-Thesis}. 

Assume that $\phi^\mfkd=1_{L}^\mfkd$. Then for $\xi>0$, the $\xi$-th Fourier coefficient $Z(\xi, \phi)$ admits natural admissible extensions (via moduli descriptions) to the RSZ integral model (here we use Bruiner's admissible Green function  $\mathcal{G}^{\mathbf B}(\xi, \phi)$ in \cite{bruinier2012regularized})
\[
\widehat{\mathcal{Z}}^{\mathbf B}(\xi, \phi) \to \mathcal{M}.
\]
Here the Cartier divisor $Z(\xi,\phi)$ on $M(V)$ is related to the Fourier expansion of $Z(h_\infty,\phi), h_\infty \in \SL_2(\mathbb R)$ in adelic language via weight $n$ Whittaker function
\[
W_\xi^{(n)}(h_\infty)=|a|^{n/2} e^{2\pi i \xi(b+ai)} e^{in\theta}, \quad \forall h_\infty=\begin{pmatrix} 1 & b \\ 0  & 0\end{pmatrix}\begin{pmatrix} a^{1/2} & 0 \\ 0  & a^{1/2} \end{pmatrix}\begin{pmatrix} \cos \theta & \sin \theta \\ -\sin \theta & \cos \theta \end{pmatrix}.
\]
We define $Z(0,\phi)=-\phi(0)c_1(\omega) \in \mathrm{Ch}^1(M(V))$, where $\omega$ is the Hodge line bundle (line bundle of weight $1$ modular forms, which is ample when the level is sufficiently small). Note that $\omega$ is also called the tautological line bundle $\CL:=(\Lie A)^\vee_{rk=1}$ (negative part of $\Lie A$ under the $F$-action of signature $(n-1,1)$) as the descent of the tautological line bundle on $\mathbb D_{n-1}$, see \cite{2017modularity}. For our purpose, we ignore the issue of definitions of integral models of $\omega$ (see e.g. \cite{howard2019linear}). For $h=(h_\infty, h_f) \in \SL_2(F_{0,\infty}) \times \SL_2(\mathbb A_f)$, we have
\begin{equation} \label{adelic to Fourier}
Z(h,\phi)=Z(0,\omega(h_f)\phi) W_0^{(n)}(h_\infty) + \sum_{\xi>0} Z(\xi,\omega(h_f)\phi) W_\xi^{(n)}(h_\infty).    
\end{equation}

\begin{definition}
For $\xi >0$, define the arithmetic intersection number (here $\tau(Z^\BQ)$ is a factor of components for the toric $Z^\BQ$ part of RSZ Shimura datum)
$$
\Int(\alpha, \xi, \phi, \varphi):= \frac{1}{\tau(Z^\BQ)}(\widehat{\mathcal{Z}}^{\mathbf B}(\xi, \phi), {}^{\mathbb L} \mathcal{TCM}_{\alpha}(\varphi))_{Ara},
$$
which is canonically decomposed into local terms (as the generic fiber has empty intersection by complex uniformization)
\begin{equation} \label{eq: deom to w-part geometric}
    \Int(\alpha, \xi, \phi, \varphi)=\Int(\alpha, \xi, \phi, \varphi)_{Ara,\infty} + \sum_{w \nmid \infty} \Int(\alpha, \xi, \phi, \varphi)_{Ara,w} \in \Int(\alpha, \xi, \phi, \varphi)_{Ara,\infty} + \sum_{\ell \not \in \mfkd} \mathbb Q \log \ell.
\end{equation}

\end{definition}

\subsection{Step II: non-archimedean intersections and uniformizations}

From moduli descriptions, we have $\Int(\alpha, \xi, \phi, \varphi)_{Ara,w}=0$ unless $w$ is inert in $F$ and $\xi \geq 0$, in which case the intersection has support in the basic locus of the mod $w$ fiber of $\CM$. Assume that $w \not \in \mfkd$ is inert in $F$. Via moduli descriptions, we now describe the basic uniformization of the embedding (which is an adelic description away from $w$)
\[
\CM \to \CM^{E_0}.
\]
After base change to $\breve{\BZ}_w$, formal completion along the mod $w$ basic locus and working on fixed components in the toric $Z^\BQ$-part, we have a closed embedding
\[
\CM^{\wedge}_0 \to \CM^{E_0, \wedge}_0.
\]
Let $V^{(w)}$ be the totally definite nearby hermitian space of $V$ at $w$. We have two embeddings $H=\U(V) \to G=\Res_{E_0/F_0}\U(V_{E_0})$ and $H^{(w)}:= \U(V^{(w)}) \to G^{(w)}=\Res_{E_0/F_0}\U(V^{(w)}_{E_0})$ which agree way from $w$. We have uniformization isomorphisms
\[
G^{(w)}(\mathbb Q) \backslash [\mathcal{N}_{n,w} \times G(\mathbb A_f^w)/K^w_G] \cong \CM^{E_0, \wedge}_0,
\]
\[
H^{(w)}(\mathbb Q) \backslash [\mathcal{N}_{L_w} \times H(\mathbb A_f^w)/K^w] \cong \CM^{\wedge}_0.
\]
Here $\mathcal{N}_{n,w}= \mathcal{N}_{L_w} \times \mathcal{N}_{L_w}$ (resp. $\CN_{n,w}=\CN_n^\GL$) if $w$ is split in $E_0$ (resp. $w$ is inert in $E_0$). Consider the discrete set 
	\[
	\Hk_{g_\mfkd}^{(w)}= \{ (g_1,g_2) \in (G(\BA_f^{w})/K^{w}_G )^2 | g_1^{-1}g_2 \in K_Gg_\mfkd K_G \}. 
	\]
We have uniformization of Hecke correspondences:
\[ 
		G^{(w)}(\mathbb Q) \backslash [ \CN_{n,w} \times \Hk_{g_\mfkd}^{(w)}]  \cong \Hk_{g_\mfkd,0}^{\wedge} \to \CM^{E_0,\wedge}_0 \times \CM^{E_0,\wedge}_0.
\]    
Hence the relative translation cycle $\pr_1: \mathcal{TCM}(g_\mfkd)^{\wedge}_0 \to \CM_0^\wedge$ over the component $\CM_0^\wedge$ is
\[
G^{(w)}(\mathbb Q) \backslash \{ (x, g_1, g_2) \in \CN_{n,w} \times \Hk_{g_\mfkd}^{(w)} | \exists (y_i, h_i) \in \CN_{n,L} \times H(\mathbb A_f^w)/K^w,  (x, g_i) \in G^{(w)}(\mathbb Q).(y_i, h_i), i=1,2  \}.
\]
Assume that $w=\ell$ is inert in $E_0$. Similar to \cite[Proposition 7.17 (2)]{Proof-AFL} (where $w=\ell$ is split in $E_0$), uniformization gives the equality of derived $1$-cycles (over formal schemes)
\begin{equation}\label{unif: CM(alpha)}
    {}^\BL \mathcal{TCM}_{\alpha}(\varphi)^{\wedge}_0= \sum_{(\delta, h) \in H^{(w)}(\mathbb Q) \backslash [B^w(\alpha)(\mathbb Q) \times H(\mathbb A^w_f) / K^w]}                             (\int_{H(\mathbb A_f^w)}  \varphi^{w}(h^{-1}\delta h_2) dh_2) [ \CN_L^{\Herm}(\delta) \times 1_{hK^{w}}].
\end{equation}
Here $B^w(\alpha)(\mathbb Q) \subseteq G^{(w)}(\mathbb Q) / H^{(w)}(\mathbb Q)$ is the preimage of $\alpha \in [H^{(w)} \backslash  G^{(w)}/ H^{(w)}](\mathbb Q)$. The orbital integral (we regard $G^{(w)}(\mathbb A_f^w) = G(\mathbb A_f^w)$)
$$
\int_{H(\mathbb A_f^w)}  \varphi^{w}(h^{-1}\delta h_2) dh_2
$$ only depends on $\delta \in G^{(w)}(\mathbb Q)/H^{(w)}(\mathbb Q)$. By moduli description of uniformization isomorphisms, we have uniformization of Kudla-Rapoport divisors when $\xi >0$ (so $u \not =0$):
\begin{equation} \label{unif: Z(xi, phi)}
    \mathcal{Z}(\xi, \phi)^\wedge_0=\sum_{(u, h') \in H^{(w)}(\BQ) \backslash [ V^{(w)}(\mathbb Q) \times H(\BA^{w}_f)/K^w], (u,u)=\xi} \phi^{w}(h'^{-1}u) [\CZ(u) \times 1_{h'K^w}] \to \CM^\wedge_0.
\end{equation}
Consider product of local non-archimedean orbital integrals away from $w$:
\begin{equation} \label{prime-to-w orb to be matched to analytic side}
\Orb^{w}((g,u), \varphi^w \otimes \phi^w) = \int_{H(\mathbb A_f^w)} \varphi^w(h^{-1}g h_2) \phi^w(h^{-1}u) dh_2 dh,   
\end{equation}

From (\ref{unif: CM(alpha)})(\ref{unif: Z(xi, phi)}), by summation over $(h,h')$ as in \cite[Theorem 9.4]{Proof-AFL} we have
\begin{equation}\label{eq: local-global decomp}
\Int(\alpha, \xi, \phi, \varphi)_{Ara,w}= 2 \log \ell \sum_{(g, u) \in H^{(w)}(\mathbb Q) \backslash B^w(\alpha) (\mathbb Q) \times V^{(w)}(\mathbb Q), (u,u)=\xi} \Height^{\Herm, \BV}_w(g, u) \Orb^{w}((g,u), \varphi^w \otimes \phi^w),    
\end{equation}

When $w=\ell$ is split in $E_0$, above formula (\ref{eq: local-global decomp}) still holds after some modifications into the AFL set up, see \cite[Theorem 9.4]{Proof-AFL}. 

As we fix regular semi-simple $\alpha$, we only need to consider orbital integrals for $(g,u) \in B^w(\alpha)(\mathbb Q) \times V^{(w)}(\mathbb Q)$. We have partial transfers for these orbital integrals at each finite place (see \cite[Definition 13.1, Remark 13.2]{Proof-AFL}), which is easily reduced to the $n=1$ case. Then we can match these orbital integrals (\ref{prime-to-w orb to be matched to analytic side}) with orbital integrals on the analytic side (in Step III).

\subsection{Step III: globalization of orbital integrals and analytic generating functions from relative trace formulas}

The goal now is to compare $\Int(\alpha, \xi_0, \phi, \varphi)_{Ara,w}$ with $w$-part of $\xi_0$-the Fourier coefficient of an analytic function $2 \partial J(\alpha, h,\phi', \varphi')$ on $h \in \SL_2(\mathbb A_{F_0})$, which is related to local derived orbital integrals at $w$ (for Conjecture \ref{TAT(g,u) conj}) and above local orbital integrals (\ref{prime-to-w orb to be matched to analytic side}) away from $w$. Here we have to transfer $\alpha \in [ \U(V^{(v_0)})  \backslash \U(V^{(v_0)}_{E_0})/ \U(V^{(v_0)}) ](F_0)$ into an invariant on the analytic side. This motivates us to consider the action of $\GL_n(F_0)$ on $\GL_n(F')/\GL_n(F_0)$ where $F'$ is the third subfield of $E$ besides $E_0$ and $F$, then we have a matching of GIT quotients and regard $\alpha \in [\GL_{n,F_0} \backslash \GL_{n,F'}/\GL_{n,F_0}](F_0)$ to define $\alpha$-part of orbits and orbital integrals.

We construct a analytic generating function $J(\alpha, h, \phi', \varphi', s), s \in \mathbb C$  similar to \cite[(11.20)(11.21)]{Proof-AFL}, via $\alpha$-sliced version of relative trace formulas (with suitable normalizations at nilpotent orbits via Tate thesis \cite[(12.20)]{Proof-AFL}). Its derivative at $s=0$ gives $\partial J(\alpha, h,\phi', \varphi')$. By Leibniz rule, we have local-global decomposition (see \cite[(12.32)(12.36)(12.38)]{Proof-AFL} \cite[\S 7.3]{Proof-Invent-cycles-TAFL}) where $\partial J_0$ is a nilpotent term:
\begin{equation}\label{eq: deom to w-part analytic}
    \partial J(\alpha, h,\phi', \varphi')=\partial J_0(\alpha, h,\phi', \varphi') + \partial J_\infty(\alpha, h,\phi', \varphi') + \sum_{w \nmid \infty} \partial J_w(\alpha, h,\phi', \varphi').
\end{equation}
\begin{equation} \label{eq: local-global decomp analytic}
 \partial J_w(\alpha, h,\phi', \varphi')=\sum_{(\gamma,u_1,u_2)} \partial \Orb((\gamma,u_1,u_2), \varphi'_w, \omega(h)\phi'_w) \Orb((\gamma,u_1,u_2), \varphi'^w, \omega(h) \phi'_w ).
\end{equation}
In particular, when $w=v_0$ we find derived orbital integrals in Conjecture \ref{TAT(g,u) conj}.

Then we plug in Gaussian
test functions at archimedean places. In other words, set $\varphi_\infty(g)=1_{\U(n)(\mathbb R \times \mathbb R)}$ and $\phi_\infty(u)=e^{-\pi (u,u)}$ for $(g,u) \in (\U(V^{(v_0)}_{E_0})\times V^{(v_0)} )(\mathbb R)=\U(n) \times \U(n) \times \mathbb R^n$. Under the Weil representation $\omega$, we have $\omega(h_\infty)\phi_\infty(u)=e^{i n \theta} |a|^{1/2} e^{\pi i (b+ia)(u,u)}$ under Iwasawa decomposition of  $h_\infty$. For fixed $\alpha$, by purely archimedean computations, we could produce partial Gaussian test function $(\phi'_\infty, \varphi'_\infty)$ on the analytic side whose orbital integrals matches $(\varphi_\infty, \phi_\infty)$ for $\alpha$-part orbits. And $(\phi'_f,\varphi'_f)$ is a partial transfer of $(\phi,\varphi)$ for $\alpha$-part orbits. In other words, test functions $(\phi',\varphi')$ are Gaussian partial transfers of $(\phi,\varphi)$, see \cite[\S 14]{Proof-Invent-cycles-TAFL}. Note that $(\phi',\varphi')$ is incoherent as its orbital integrals transfer to orbital integrals for the incoherent hermitian space $V^{(v_0)}_\mathbb R \otimes V(\mathbb A_f)$.

Via (\ref{adelic to Fourier})(\ref{eq: local-global decomp})(\ref{eq: local-global decomp analytic}) and (\ref{eq: deom to w-part geometric})(\ref{eq: deom to w-part analytic}), we need to do local-global comparison of $\Int(\alpha, \xi, \phi, \varphi)_{Ara,w}$ with $-2 \partial J_w(\alpha, \xi, \phi', \varphi')$ for $\xi>0$. Firstly we need an archimedean comparison (Step IV).

\subsection{Step IV: archimedean intersections, archimedean orbital integrals and holomorphic differences}

In the case $n=1$, we can compute archimedean orbital integrals for Gaussian test function $\phi'_\infty(u_1,u_2)=2^{-3/2}(u_1+u_2)e^{- \frac{1}{2} \pi(u_1^2+u_2^2)}$ directly in terms of K-Bessel function $K_s(c)=\frac{1}{2}\int_{\mathbb R_{>0}} e^{\frac{-1}{2} c (u+1/u)} u^s \frac{du}{u}, c>0, s \in \mathbb C$, which shows a potential relation to Kudla's Green functions after taking derivative at $s=0$, see \cite[Lemma 12.3]{Proof-AFL}. For general $n$, we use \cite[(12.11))]{Proof-AFL} to reduce to $n=1$. By such purely archimedean direct computations \cite[Corollary 10.3, Lemma 14.3-14.4]{Proof-AFL}, the analytic archimedean term $\partial J_\infty(\alpha, h,\phi', \varphi')$ is related to archimedean intersection numbers, i.e. values of Kudla's Green functions \cite[(8.13)]{Proof-AFL} on the complex fiber $M(V)(\mathbb C)$:
$$
\CG^{\bf K}(\xi, h_\infty, \phi)=\sum_{(g,u)} \phi(g^{-1}u) (\CG^{\bf K}(u,h_\infty) \times 1_{gK}).
$$ 
Here $\CG^{\bf K}(0,h_\infty)=- \log |a|$ for $h_\infty=\begin{pmatrix} 1 & b \\ 0  & 0\end{pmatrix}\begin{pmatrix} a^{1/2} & 0 \\ 0  & a^{1/2} \end{pmatrix}\begin{pmatrix} \cos \theta & \sin \theta \\ -\sin \theta & \cos \theta \end{pmatrix}$. By complex uniformization, the geometric archimedean term
\begin{equation}
 \Int(\alpha, \xi, \phi, \varphi)_{Ara,\infty}=\sum_{(g,u) \in H(\mathbb Q) \backslash B(\alpha)(\mathbb Q) \times V(\mathbb Q),  (u,u)=\xi} \Int_\infty(g,u) \Orb((g,u),\varphi \otimes \phi)   
\end{equation}
where $\Int_\infty(g,u)=\CG^{\bf K}(u,h_\infty)(z_g)$ where $z_g \in \mathbb D_{n-1}$ is the unique fixed point of $g \in B(\alpha)(\mathbb R) \subseteq \U(V_\infty)$. The work \cite{ES-differenceGreen} (see also \cite{Buck2023automorphic} for Hirzebruch–Zagier divisors) shows that the generating function of differences of Green functions 
$$
\widehat{\mathcal{Z}}^{\mathbf K -\mathbf B}(h,\phi)= \sum_{\xi \in \mathbb Q} ( \CG^{\bf K}(\xi, \omega(h_f)\phi, h_\infty) -\CG^{\bf B}(\xi, \omega(h_f)\phi) ) W^{(n)}_\xi(h_\infty) 
$$ is a weight $n$ modular form on $h=(h_\infty, h_f) \in \SL_{2}(\mathbb A_{F_0})$. Set $\Int^{\mathbf K - \mathbf B} (\alpha, h, \phi, \varphi)= \frac{1}{\tau(Z^\BQ)}\widehat{\mathcal{Z}}^{\mathbf K -\mathbf B}(h,\phi)(\mathrm{TCM}_\alpha(\varphi)(\mathbb C))$. Via purely archimedean computation as in \cite[Proposition 14.5.]{Proof-AFL}, we find $2 \partial J_\infty(\alpha, \xi,\phi', \varphi') = - \Int^{\mathbf K}_\infty (\alpha, \xi, \phi, \varphi)$ for $\xi>0$, and
\[
 2 \partial J(\alpha, h,\phi', \varphi') + \Int^{\mathbf K - \mathbf B} (\alpha, h, \phi, \varphi)
\]
is a holomorphic modular form of weight $n$ on $h \in \SL_2(\mathbb A_F)$.

\subsection{Step V: modularity of geometric theta series and (holomorphic) modification of $1$-cycles}

To prove arithmetic transfers over general $p$-adic fields, we do not know modularity of arithmetic theta series over integral models with corresponding parahoric levels. Instead, as in \cite[\S 12-13]{Proof-AT} we find another $1$-cycle (modified twisted derived CM cycles)
\[
{}^{\mathbb L} \mathcal{TCM}_{\alpha}(\varphi)^{mod}={}^{\mathbb L} \mathcal{TCM}_{\alpha}(\varphi)+\mathcal{C}
\]
with some explicit and computable (``very special'') $1$-cycle $\mathcal C$ such that the pairing
\[
(-, {}^{\mathbb L} \mathcal{TCM}_{\alpha}(\varphi)^{mod})_{Ara}
\]
factors through the natural map $\Ch^{1,adm}(\CM) \to \mathrm{Ch}^1(M(V))$. Then $(Z(h,\phi), {}^{\mathbb L} \mathcal{TCM}_{\alpha}(\varphi)^{mod})_{Ara})$ is a well-defined holomorphic modular form by modularity over reflex field \cite{Liu-Thesis}. Define the modified difference 
\[
\mathcal{E}^{mod}(h)= 2 \partial J(\alpha, h,\phi', \varphi')^{mod} + \Int^{\mathbf K - \mathbf B} (\alpha, h, \phi, \varphi)^{mod} + (Z(h,\phi), {}^{\mathbb L} \mathcal{TCM}_{\alpha}(\varphi)^{mod})_{Ara},
\]
by subtracting contributions from $\mathcal C$, which is now a holomorphic modular forms of weight $n$ on $h \in \SL_2(\mathbb A_{F_0})$. If Conjecture \ref{TAT(g,u) conj} holds, then $\mathcal{E}^{mod}(h)$ is a constant by local-global relations. To prove Conjecture \ref{TAT(g,u) conj}, the next step is to show that (the $v_0$-part of) $\mathcal{E}^{mod}(h)$ is indeed a constant.

\subsection{Step VI: matching totally positive Fourier coefficients by induction, simple known cases and modularity}

We use two facts.

\begin{enumerate}
    \item We crucially use mirabolic special cycles to show that TAFL in rank $n-1$ (which we assume by induction) implies TAFL$(g,u)$ when $(u,u)$ is a unit. See \cite[\S 5.3-5.4]{Proof-Invent-cycles-TAFL} for arithmetic induction on geometric and analytic sides. Hence by induction and maximal order cases of TAFLs, we find that $\xi$-th Fourier coefficient of $\mathcal{E}^{mod}(h)$ is zero for $\xi>0$ with $v_p(\xi)=0$.
    \item Above modularity results also explicitly describe the levels of holomorphic modular forms, which implies that $\mathcal{E}^{mod}(h)$ has maximal level at $p$ (resp. Iwahori level at $p$) for the proof of TAFL (resp. TAT).
\end{enumerate}

 In the TAFL case, from knowledge of modular forms  \cite[Lemma 13.6]{Proof-AFL}, these two facts imply that $\mathcal{E}^{mod}(h)$ is a constant. In the TAT case, the modification over mod $p$ fibers and complex fibers is the same as \cite[\S 12-13]{Proof-AT}, and $\mathcal{E}^{mod}(h)$ is a holomorphic modular form (thanks to holomorphic modularity of specific theta series \cite[Theorem 13.9]{Proof-AT}). We apply the double induction method \cite{Proof-AT} and work with similar difference functions defined via geometric theta series for the Fourier transform of $\phi$ at $v_0$, to see that $\mathcal{E}^{mod}(h)$ is a constant.

\subsection{Step VII: shrinking supports and from global identities to local TAFLs}

Note after subtracting the archimedean terms,
the intersection numbers and derived orbital integrals at non-archimedean places all lie in $\mathbb Q$-linear span of log $\ell$ ($\ell \not \in \Delta$). Using $\mathbb Q$-linearly independence of $\log \ell$, from the vanishing of $\xi_0$-th Fourier coefficients of  $\mathcal{E}^{mod}(h)$ we obtain a semi-global identity (adding these modification terms back)
$$ 
2 \partial J(\alpha, \xi_0,\phi', \varphi')_{v_0} + \Int^{\mathbf K - \mathbf B} (\alpha, \xi_0, \phi, \varphi)_{v_0} + \Int(\alpha, \xi_0, \phi, \varphi)_{Ara,v_0} =0.
$$
By shrinking supports at another split place of $F/F_0$ (applying the technique of \cite[Lemma 13.7]{Proof-AFL}), we may assume the sum for $\Int(\alpha, \xi_0, \phi, \varphi)_{Ara,v_0}$ has only one non-zero term (and is the term for orbits of $(g_0,u_0)$) and the orbital integral $\Orb^{v_0}(g_0,u_0), \varphi^{v_0} \otimes \phi^{v_0})$ is non-zero. Then we obtain the TAFL for $(g_0, u_0)$ from above semi-global identity \cite[\S 8]{Proof-Invent-cycles-TAFL}. The TAT conjecture for vertex lattices may be proved similarly.

\bibliographystyle{alpha}
\bibliography{reference}

\end{document}